\documentclass[11pt]{article}

\usepackage[T1]{fontenc}
\usepackage[utf8]{inputenc}
\usepackage{lmodern}
\usepackage[a4paper,margin=1in]{geometry}
\usepackage{amsmath,amssymb,amsthm,mathtools}
\usepackage{aliascnt}
\usepackage{enumitem}
\usepackage[protrusion=true,expansion=false]{microtype}
\usepackage[dvipsnames]{xcolor}
\usepackage[colorlinks=true,linkcolor=MidnightBlue,citecolor=MidnightBlue,
  urlcolor=MidnightBlue]{hyperref}

\makeatletter
\newcommand\blfootnote[1]{%
  \begingroup
  \renewcommand\@makefntext[1]{\noindent\strut ##1}%
  \renewcommand\thefootnote{}\footnote{#1}%
  \addtocounter{footnote}{-1}%
  \endgroup
}
\makeatother

\setlist{itemsep=0.25em,topsep=0.4em}
\newtheorem{theorem}{Theorem}[section]
\newaliascnt{proposition}{theorem}
\newtheorem{proposition}[proposition]{Proposition}
\aliascntresetthe{proposition}
\newaliascnt{lemma}{theorem}
\newtheorem{lemma}[lemma]{Lemma}
\aliascntresetthe{lemma}

\newaliascnt{corollary}{theorem}

\aliascntresetthe{corollary}
\theoremstyle{definition}
\newaliascnt{definition}{theorem}

\aliascntresetthe{definition}

\newaliascnt{remark}{theorem}

\aliascntresetthe{remark}

\newaliascnt{example}{theorem}

\aliascntresetthe{example}

\usepackage[nameinlink,capitalise,noabbrev]{cleveref}

\newcommand{\Rm}{\operatorname{Rm}}
\newcommand{\Ric}{\operatorname{Ric}}
\newcommand{\Fix}{\operatorname{Fix}}
\newcommand{\Vol}{\operatorname{Vol}}
\newcommand{\Length}{\operatorname{Length}}
\newcommand{\Isom}{\operatorname{Isom}}
\newcommand{\Diff}{\operatorname{Diff}}
\newcommand{\AVR}{\operatorname{AVR}}

\newcommand{\cR}{\mathcal R}
\newcommand{\R}{\mathbb R}
\newcommand{\Sph}{\mathbb S}
\newcommand{\id}{\mathrm{id}}

\newcommand{\Scal}{\operatorname{Scal}}
\newcommand{\Hess}{\operatorname{Hess}}
\newcommand{\Fl}{\operatorname{Fl}}

\title{\textbf{Fixed geodesics of involutions on $3$-spheres}}
\author{Bin Shen and Hui Zhang}
\date{}

\begin{document}

\maketitle

\begin{abstract}
Let an isometric involution act on a smooth Riemannian $3$-sphere with
sectional curvature at least one. We prove that every one-dimensional
connected component of its fixed-point set has length at most $2\pi$,
and that the bound is sharp. The proof realizes the fixed circle as the
link of a complete totally geodesic surface in a positively curved
expanding Ricci soliton and applies Bishop--Gromov comparison to its area
growth. The analytic input is an equivariant version of Deruelle's
smoothing of positively curved cones. We use Gianniotis--Schulze's
equivariant local deformation and explain a marked compactness argument
for arbitrary finite isotropy. In particular, the end chart retains the
specified group action, rather than only an unspecified conjugate.
\end{abstract}

\blfootnote{%
}

%\tableofcontents

\section{Introduction}

Petrunin asks whether a closed geodesic that is the fixed-point set of an
isometric involution of a Riemannian $3$-sphere with sectional curvature
at least one has length at most $2\pi$ (cf. \cite[Section~7]{Petrunin}).
We answer the smooth Riemannian question affirmatively.

\begin{theorem}\label{thm:main}
Let $(\Sph^3,g)$ be a smooth Riemannian three-sphere with $\sec_g\geq1$,
and let $\tau$ be an isometric involution. Every one-dimensional
connected component $Y$ of the fixed-point set of $\tau$, denoted by $\Fix(\tau)$, satisfies
\[
 \Length_g(Y)\leq2\pi.
\]
The bound is attained by a great circle fixed by a linear involution of
the unit round sphere.
\end{theorem}

We specify our curvature convention before stating the general result.
For any Riemannian metric $h$, the induced metric on $\Lambda^2T_xX$ is determined by
\[
 \langle u\wedge v,w\wedge z\rangle_h
 =h(u,w)h(v,z)-h(u,z)h(v,w).
\]
The \emph{curvature operator} is the self-adjoint endomorphism
$\cR_h:\Lambda^2T_xX\to\Lambda^2T_xX$ defined by
\begin{equation}\label{eq:curvature-definition}
 \langle\cR_h(u\wedge v),w\wedge z\rangle_h
 =h(R^h(u,v)z,w),
\end{equation}
where $R^h$ is the curvature tensor of $h$. 
We use the curvature-operator convention of Petersen~\cite[Sections~3.1.2--3.1.3]{Petersen}, normalized so that the unit round sphere has $\mathcal R_h=I$, where
$I=\id_{\Lambda^2TX}$. The inequality $\cR_h\geq I$ means
\[
 \langle\cR_h\omega,\omega\rangle_h\geq|\omega|_h^2
 \quad\text{for every }x\in X,\ \omega\in\Lambda^2T_xX.
\]
It is well known that
the sectional curvature is
$\langle\cR_h(u\wedge v),u\wedge v\rangle_h$, for any orthonormal $u,v$.
Sectional curvature therefore only tests decomposable bivectors, whereas the inequality $\cR_h\geq I$ tests all bivectors. However, in dimension three every
bivector is decomposable, so $\sec_h\geq1$ is equivalent to
$\cR_h\geq I$. This equivalence is not asserted in higher dimensions.
We denote the unit round metric by  $g_{\rm rd}$.

\begin{theorem}\label{thm:fixed-link}
Let $(X^m,h)$ be closed and simply connected, with $m\geq2$ and
$\cR_h\geq I$, and let a finite group $G$ act isometrically on $X$.
If $Y^d$ is a connected component of the common fixed set $X^G$ and
$d\geq1$, then
\[
 \Vol_h(Y)\leq\Vol(\Sph^d,g_{\rm rd}).
\]
\end{theorem}

For $d\geq2$, since $Y$ is totally geodesic, the Bishop
comparison already gives this bound. The essential application is $d=1$.
A compact connected one-dimensional fixed component is a circle.
%Its cone has one radial variable and one angular variable, hence dimension two. 
If the circle has length $L$, the area of its cone up to radius $R$
is $LR^2/2$. The desired length bound will consequently follow from the
upper area-growth bound $\pi R^2$ on a complete nonnegatively curved
surface.

To make Bishop--Gromov available, we smooth the cone over $X$ by an equivariant expanding gradient Ricci soliton. The smoothing used to produce that surface must respect the specified action on the link, including its fixed points.

\begin{theorem}[Equivariant smoothing]\label{thm:equivariant-smoothing}
Let $(X^m,h)$ be closed and simply connected, where $m\geq2$, and suppose
$\cR_h>I$. Let $\rho:G\to\Isom(X,h)$ be an action of a finite group.
There are a complete expanding gradient Ricci soliton
$(M^{m+1},g,\nabla f)$ with positive curvature operator, an isometric
action $\widehat\rho:G\to\Isom(M,g)$, and an asymptotic chart
\[
 \Phi:(R_0,\infty)\times X\longrightarrow M\setminus K,
\]
for some compact set $K$, such that after translating $f$ by a constant,
\begin{equation}\label{chart}
	f\circ\Phi=\frac{r^2}{4},\qquad
	\widehat\rho(\sigma)\circ\Phi
	=\Phi\circ(\id\times\rho(\sigma))\quad(\sigma\in G).
\end{equation}
The chart %has the rate-two estimates of \eqref{eq:rate-two} and 
can be
chosen in the Morse-flow gauge
\begin{equation}\label{eq:intro-morse}
 d\Phi(\partial_r)
 =\frac r2\frac{\nabla f}{|\nabla f|^2}\circ\Phi.
\end{equation}
In particular, $f$ is $G$-invariant and the marking at infinity is the
original action $\rho$.
\end{theorem}

Deruelle proved the nonequivariant smoothing theorem in
\cite[Theorem~1.3]{DeruelleGAFA}. Gianniotis--Schulze established an
equivariant version for a finite orthogonal action free on the spherical
link \cite[Section~6]{GianniotisSchulze}. Their openness construction uses
equivariant normal graphs and flows and does not require freeness (cf. Subsection \ref{subsection 3.1}).
Their closedness argument, however, uses quotient orbifolds with an
isolated singularity. It does not directly cover an action with fixed
points on the link, because an element fixing $x\in X$ fixes the whole ray
$\Phi((R_0,\infty)\times\{x\})$, so the quotient has singular strata
extending to infinity. Thus we cannot simply cite that theorem for
\cref{thm:equivariant-smoothing}.

We supply the required marked closedness on the smooth covering
manifolds. Its central point is to combine the convergence on fixed finite
annuli with the \emph{uniform} decay of the rescaled level metrics. This
identifies the limiting link marking. Analytic continuation then extends
the limiting action from the end to the whole expander.

The remainder of the manuscript is organized as follows. \Cref{sec:inputs} collects the definitions and cited analytic
results. \Cref{sec:smoothing} proves equivariant smoothing, with expanded
openness and closedness arguments. \Cref{sec:fixed} then proves
\cref{thm:fixed-link,thm:main} by the geometric comparison.

\section{Preliminaries}\label{sec:inputs}

\subsection{Cones, asymptotic estimates, and fixed cylinders}

For a closed Riemannian manifold $(X^m,h)$, the punctured metric cone is
\[
 C(X)=(0,\infty)\times X,\qquad g_C=dr^2+r^2h.
\]
The metric completion adds the cone vertex. Under the natural
orthonormal identification of tangential bivectors, the warped-product
curvature formula is
\begin{equation}\label{eq:cone-curvature}
 \cR_{g_C}|_{\Lambda^2TX}=r^{-2}(\cR_h-I),\qquad
 \cR_{g_C}(\partial_r\wedge\cdot)=0.
\end{equation}
Thus $\cR_h>I$ gives a nonflat cone with nonnegative curvature operator.
Constant scaling gives $\cR_{c^2h}=c^{-2}\cR_h$ as endomorphisms.

A \emph{marked asymptotic chart} of a complete manifold is a
diffeomorphism
\[
 \Phi:(R_0,\infty)\times X\longrightarrow M\setminus K
\]
for a compact $K$, together with the specified link $(X,h)$.
We use the usual exterior-chart convention, i.e., for every sufficiently large
$R_1$, the complement of $\Phi((R_1,\infty)\times X)$ is compact.
Deruelle's sharp asymptotic estimates \cite[Theorem~1.2 and proof of Theorem~3.2]{DeruellePisa} gives the rate-two estimates of the chart in \eqref{chart}, that is
\begin{equation}\label{eq:rate-two}
 \left|(\nabla^{g_C})^j e\right|_{g_C}\leq C_jr^{-2-j},
 \qquad e=\Phi^*g-g_C,\qquad j=0,1,2,\ldots.
\end{equation}
They imply smooth conical asymptotics. The stronger weighted decay in
\eqref{eq:rate-two}, not only unweighted convergence to zero, will be used
in the analytic argument (cf. Subsection \ref{subsecction2.4}). 

The $g_C$ in \eqref{eq:rate-two} is called the \emph{relative norm}, which is induced by the cone metric. In local coordinates,
\[
 |e|_{g_C}^2=(g_C)^{\alpha\gamma}(g_C)^{\beta\delta}e_{\alpha\beta}e_{\gamma\delta}.
\]
In particular, $|e|_{g_C}\leq\varepsilon$ implies
\begin{equation}\label{eq:metric-sandwich}
 (1-\varepsilon)g_C\leq\Phi^*g\leq(1+\varepsilon)g_C.
\end{equation}
For $j\geq 1$, \eqref{eq:rate-two} is also named as the \emph{differentiated estimates}, since it takes iterated covariant derivatives with the cone connection. 

If $Y\subset X$ is a connected fixed component, the \emph{displayed cylinder} of $Y$ is 
\[
 E_Y=\Phi((R_1,\infty)\times Y),\qquad R_1>R_0.
\]
It is a cylinder in the topological sense, namely, a
product of a radial interval with $Y$. But we adopt a conical metric
$dr^2+r^2h|_Y$ as its model metric, not the product metric $dr^2+h|_Y$.
Since $Y$ is totally geodesic in $X$ (see \cref{lem:fixed-submanifold}),
$C(Y)$ is totally geodesic in $C(X)$. Restricting $e$ and its tangential
covariant derivatives to $C(Y)$ therefore preserves
\eqref{eq:rate-two}, relative to $dr^2+r^2h|_Y$.

\subsection{Asymptotic volume ratio}\label{subsec:avr}

For a complete connected $n$-manifold $(P,k)$ with $\Ric_k\geq0$, define
\begin{equation}\label{eq:avr}
 \AVR(k)=\AVR(P,k)
 :=\lim_{R\to\infty}\frac{\Vol_k B_k(p,R)}{R^n}.
\end{equation}
The limit exists because the sequence is nonincreasing by Bishop–Gromov comparison.
Noticing there is no factor
$\omega_n^{-1}=\Vol^{-1}(B_{\R^n}(0,1))$ in \eqref{eq:avr}, one could deduce that
\[
 0\leq\AVR(k)\leq\omega_n,\qquad \AVR(g_{\R^n})=\omega_n.
\]
It is independent of $p$. Indeed, for $D=d_k(p,q)$ and $R>D$, we have
\[
 B_k(p,R-D)\subset B_k(q,R)\subset B_k(p,R+D).
\]
Division by $R^n$ and passage to the limit prove the assertion.

For a cone of dimension $n$, the volume element and annular volume are
\[
 d\mu_{g_C}=r^{n-1}\,dr\,d\mu_h,\qquad
 \Vol_{g_C}\{R_1<r<R\}
 =\frac{\Vol_h(X)}{n}(R^n-R_1^n).
\]
Consequently, for an asymptotically conical end exhausting the complement
of a compact set,
\begin{equation}\label{eq:avr-link}
 n\AVR(k)=\Vol_h(X).
\end{equation}
To justify using geodesic balls here, apply
\eqref{eq:metric-sandwich} on a sufficiently far tail. Radial paths give
$d_k(p,\Phi(r,x))\leq C+\sqrt{1+\varepsilon}\,r$ uniformly in $x$.
Conversely, any path from the compact interior to level $r$ has length
at least $\sqrt{1-\varepsilon}(r-R_1)$ on its tail.
The corresponding annular volume bounds imply \eqref{eq:avr-link}
after $\varepsilon\downarrow0$.
If only one of several ends is specified, this argument gives the lower
bound $n\AVR(k)\geq\Vol_h(X)$ for that end, which is all that the geometric
proof needs (cf. the proof of Proposition \ref{prop:marked-closedness}).

\subsection{Expander potentials and Morse-flow coordinates}

We use the expander convention
\begin{equation}\label{eq:soliton}
 \Ric_g+\frac12g=\Hess_g f.
\end{equation}
Adding a constant to $f$ changes neither this equation nor $\nabla f$.
For compactness, choose the integral normalization
\begin{equation}\label{eq:integral-normalization}
 \int_M e^{-f}\,d\mu_g=(4\pi)^{n/2}.
\end{equation}
The Gaussian integral is finite for the conical expanders in question.
This normalization removes arbitrary additive constants along a sequence.

There is a different normalization useful for the radial coordinate. According to the second Bianchi identity, 
the soliton identities give
\begin{equation}\label{eq:potential-constant}
 |\nabla f|^2+\Scal_g-f=c,
\end{equation}
where $c$ is constant on a connected soliton.
If $p$ is a critical point of $f$, then $c=\Scal_g(p)-f(p)$. Put
\begin{equation}\label{eq:canonical-F}
 F=f+c.
 \qquad\text{Then}\qquad
 |\nabla F|^2+\Scal_g=F.
\end{equation}
%In Deruelle's notation $c=\mu(g)$ for the integral-normalized $f$.
%The frequently used potential $v$ in \cite{DeruellePisa} is $F+n/2$, not $F$.

This explains why we later introduce such a sequence $F_i=f_i+c_i$.
The constants $c_i$ remove the constants in the scalar-gradient identities
and make
\begin{equation}\label{eq:canonical-radius}
 r=2\sqrt{F_i},\qquad
 |\nabla F_i|^2=\frac{r^2}{4}-\Scal_{g_i}
\end{equation}
hold exactly on each end (cf. Subsection \ref{subsection 3.2}). If $f_i$ were replaced by $f_i+a_i$, then
$c_i$ would be replaced by $c_i-a_i$, leaving $F_i$ unchanged.
Thus $F_i$ is canonical for this purpose. A further change $F_i\mapsto
F_i+a$ would change the squared radius by $4a$, although the gradient
and Hessian would be unchanged. 
%The two normalizations \eqref{eq:integral-normalization} and \eqref{eq:canonical-F} must therefore not be conflated.

When $\Ric_g\geq0$, \eqref{eq:soliton} gives
$\Hess F\geq g/2$. The complete conical solitons under consideration
have a proper potential and a unique nondegenerate minimum $p$.
Integration along a minimizing unit-speed geodesic from $p$ gives
$F(x)\geq F(p)+d(p,x)^2/4$. Strict convexity excludes a second critical
point. A small sublevel is a ball by the Morse lemma, and there are no
further critical points, so the gradient flow identifies all larger
sublevels with balls. In particular,
\begin{equation}\label{eq:expander-topology}
 M\cong\R^n.
\end{equation}

On a region without critical points put
\[
 W=\frac{\nabla F}{|\nabla F|^2},
\]
so that $dF(W)=1$. If $j:X\to\{F=R^2/4\}$ is a parametrization, define
\begin{equation}\label{eq:morse-formula}
 \Phi(r,x)=\Fl_W^{(r^2-R^2)/4}(j(x)),
\end{equation}
where $\operatorname{Fl}_W^t(q)$ denotes the time-$t$ flow of
the vector field $W$ starting at $q$. Since $dF(W)=1$,
the time $(r^2-R^2)/4$ carries $j(x)$ from the level
$F=R^2/4$ to the level $F=r^2/4$.
Then it satisfies
\[
 F(\Phi(r,x))=F(j(x))+\frac{r^2-R^2}{4}=\frac{r^2}{4},
\]
and the \emph{Morse-flow equation}
\begin{equation}\label{eq:morse-gauge}
 d\Phi(\partial_r)=\frac r2W\circ\Phi.
\end{equation}
The flows exist between finite regular levels by properness. Flowing
through all outward levels gives an end chart. If an isometric action
preserves $F$, then it preserves $W$ and commutes with its flow.
Consequently \eqref{eq:morse-formula} is equivariant whenever $j$ is.

The condition $F\circ\Phi=r^2/4$ alone does not imply
\eqref{eq:morse-gauge}, since radius-dependent angular reparametrizations
preserve the first condition. Conversely, if two charts use the same $F$,
both satisfy $F\circ\Phi=r^2/4$ and \eqref{eq:morse-gauge}, their transition is
\begin{equation}\label{eq:morse-transition}
 (r,x)\longmapsto(r,\psi(x))
\end{equation}
for a single diffeomorphism of cross-sections. Indeed the transition first
has the form $(r,A(r,x))$ because it preserves $F$, and differentiating the
two flow equations gives $\partial_rA=0$.
An arbitrary choice of $j$ in \eqref{eq:morse-formula} need not give a
\emph{specified} link metric in its angular coordinates. We address this
marking issue in \cref{prop:open,prop:marked-closedness}.

\subsection{Compactness, relabeling, and analytic continuation}\label{subsecction2.4}

We use the following consequence of Deruelle's compactness results. %,
%without including any assertion about a prescribed group marking.

\begin{theorem}[Deruelle compactness]\label{thm:deruelle-compactness}
Let $(M_i^n,g_i,\nabla f_i,p_i)$ ($n\geq 3$) be complete
integral-normalized asymptotically conical expanding gradient Ricci
solitons with nonnegative curvature operator, based at the unique
critical points of their potentials. Suppose their links $h_i$ on a
fixed closed manifold belong to a compact smooth subset of
$\{h:\cR_h>I\}$.
After passage to a subsequence, the pointed solitons converge smoothly
to a complete nonnegatively curved expanding gradient Ricci soliton,
and the convergence is smooth conical convergence.
The global curvature estimates used in this convergence are uniform, namely, 
on sufficiently large canonical potential levels,
\begin{equation}\label{eq:uniform-curvature}
 |\nabla^j\Rm_{g_i}|_{g_i}\leq C_jr^{-2-j}
 \quad(j\geq0),
\end{equation}
with a common lower radius and constants.
\end{theorem}

The compactness assertion is \cite[Theorem~1.7]{DeruellePisa}.
The link bounds give the required asymptotic curvature-derivative bounds
and noncollapsing, see also \cite[Remark~4.11]{DeruelleGAFA}.
Bonnet--Myers gives a uniform diameter bound for the links, and
\eqref{eq:avr-link} gives their positive asymptotic-volume-ratio bound.
The uniform global derivative estimates are supplied by
\cite[Proposition~4.3]{DeruellePisa}.
For the associated uniform tail estimates, see pp.~508--509 of \cite[the proof of Theorem~3.2]{DeruellePisa}, and equation~(4.5), p.~527 of \cite[the proof of Theorem~4.6]{DeruellePisa}.

The pointed convergence identifies compact subsets of the interiors, while the
conical convergence also involves identifications of the ends.
Compatibility of those identifications with the \emph{prescribed}
Morse-flow charts is not automatic. The fixed-level construction in
\cref{prop:marked-closedness} provides it.

\begin{lemma}[Compactness of relabelings]\label{lem:relabel}
Let $a_i\to a$ and $b_i\to b$ smoothly be metrics on closed manifolds
$X_0$ and $X$, respectively. Suppose
$\psi_i:X_0\to X$ are diffeomorphisms with $\psi_i^*b_i=a_i$.
A subsequence converges smoothly to a diffeomorphism $\psi_\infty$
satisfying $\psi_\infty^*b=a$. The inverses converge smoothly as well.
\end{lemma}

\begin{proof}
Uniform equivalence to fixed background metrics gives bounds on
$d\psi_i$ and $d\psi_i^{-1}$. In local coordinates the isometry equation
implies
\[
 \partial_\alpha\partial_\beta\psi_i^\gamma
 =\Gamma(a_i)^\delta_{\alpha\beta}\partial_\delta\psi_i^\gamma
 -\Gamma(b_i)^\gamma_{\mu\nu}(\psi_i)
       \partial_\alpha\psi_i^\mu\partial_\beta\psi_i^\nu.
\]
It bounds second derivatives, and differentiating it bounds all higher
derivatives. The first-derivative bounds permit the use of uniformly
small domain charts mapping into target charts, and the compactness supplies
finite covers. Apply the same argument to the inverses and take a
diagonal Arzel\`a--Ascoli subsequence. The limiting maps are inverses by
passage to the limit in their compositions, and the metric identity
passes to the limit as well.
\end{proof}

\begin{lemma}[Extension of end isometries]\label{lem:globalize-action}
Let $(M,g)$ be complete, connected, and real analytic, and let
$E\subset M$ be a nonempty connected open set. If
$\pi_1(E)\to\pi_1(M)$ is surjective, every isometry from $(E,g|_E)$
onto itself extends uniquely to a global isometry of $(M,g)$.
In particular this holds when $M$ is simply connected.
\end{lemma}

This is the analytic continuation theorem in 
\cite[Theorem~5.2]{KotschwarWang}. The uniqueness also follows from the fact that
isometries of a connected Riemannian manifold agreeing on a nonempty
open set agree everywhere.

Ivey's analyticity theorem
\cite[Theorem~1]{Ivey} shows that a smooth Ricci soliton metric is real analytic
in harmonic coordinates. In our sign convention its hypothesis is
\[
 \mathcal L_{\nabla F}g=2\Ric_g+g,
\]
which is equivalent to \eqref{eq:soliton}. We will use this fact in the proof of Proposition \ref{prop:marked-closedness}. %Real analyticity means that there is a compatible analytic atlas in which the metric coefficients are locally represented by convergent power series. This is a property deduced from the soliton equation, not a consequence of smooth convergence of analytic metrics.

\section{Equivariant smoothing}\label{sec:smoothing}

We prove Theorem \ref{thm:equivariant-smoothing} in this section.

\subsection{Equivariant openness}\label{subsection 3.1}

\begin{proposition}\label{prop:open}
Let $(M,g_0,\nabla F_0)$ be a positive-curvature-operator expander with
marked link $(X,\gamma_0)$ and a rate-two Morse-flow chart $\Phi_0$.
Suppose $\rho:G\to\Isom(X,\gamma_0)$ is a finite-group action extending
to $\widehat\rho_0:G\to\Isom(M,g_0)$, and $\Phi_0$ intertwines the
actions exactly. Every sufficiently small smooth $G$-invariant
deformation $\gamma$ of $\gamma_0$ with $\cR_\gamma>I$ is the marked
link of a positive-curvature-operator expander on $M$ invariant under
the same diffeomorphisms $\widehat\rho_0(\sigma)$.
Its potential and a rate-two Morse-flow chart can be chosen invariant
and exactly equivariant for the original $\rho$.
\end{proposition}

\begin{proof}
We separate the metric construction, gradient structure, and marked
coordinate construction.

\emph{The fixed action preserves the reference data.}
Write $S_\sigma=\widehat\rho_0(\sigma)$. On the end, equivariance and
$F_0\circ\Phi_0=r^2/4$ give $F_0\circ S_\sigma=F_0$.
The soliton equation and the isometry property give
\[
 \Hess_{g_0}(F_0\circ S_\sigma-F_0)=0.
\]
Then $\nabla_{g_0}(F_0\circ S_\sigma-F_0)$ is parallel and vanishes on the end,
hence vanishes everywhere. %The difference is therefore the zero constant.
Thus the reference potential and its weighted function spaces are
$G$-invariant.

\emph{Equivariance of the local inverse.}
Deruelle's deformation theorem \cite[Theorem~3.6]{DeruelleGAFA} uses a
reference extension of the link perturbation of the form
\[
 T(\gamma)=g_0+\chi(F_0)\mathcal E_{\Phi_0}(\gamma-\gamma_0),
\]
where $\chi$ is a cutoff supported on the end and equal to one farther
out, and $\mathcal E_{\Phi_0}$ is the tensor whose expression in $\Phi_0$
coordinates is $r^2(\gamma-\gamma_0)$ on angular vectors and zero on
radial vectors. In particular,
\[
 T(\rho(\sigma)^*\gamma)=S_\sigma^*T(\gamma).
\]
The Ricci tensor and the background-dependent gauge operator are natural
under pullback of all their data. Since $S_\sigma$ preserves $g_0,F_0$
and the reference extension, the gauge-fixed nonlinear map used by
Deruelle is equivariant under $S_\sigma^*$.

Write $k$ for the correction and $\mathcal Q(\gamma,k)=0$ for its
equation. Equivariance says
\[
 \mathcal Q(\rho(\sigma)^*\gamma,S_\sigma^*k)
 =S_\sigma^*\mathcal Q(\gamma,k).
\]
If $\gamma$ is invariant, $k$ and $S_\sigma^*k$ therefore solve the
same equation. Choose the local inverse neighborhood invariant by
intersecting its finitely many translates.
Uniqueness \emph{inside this gauge-fixed neighborhood} gives
$S_\sigma^*k=k$. Thus the resulting metric $g_\gamma$ is invariant
under the original diffeomorphisms $S_\sigma$ and has the specified
cone tensor $dr^2+r^2\gamma$ in the reference marking.
The positive-curvature-operator assertion follows from Deruelle's
deformation theorem with weight $\alpha=1$.
This is the equivariant deformation argument in
\cite[proof of Theorem~6.1, openness, pp.~3970--3973]{GianniotisSchulze}, and
no freeness of the link action is used in this step.
%We have used local uniqueness of a gauge-fixed correction, not unmarked global uniqueness to select a prescribed boundary isometry.

\emph{An invariant gradient potential.}
The reference expander has $\Hess F_0>g_0/2$ and is diffeomorphic to
$\R^{m+1}$ by \eqref{eq:expander-topology}.
Deruelle's gradient-structure result
\cite[Theorem~4.3]{DeruelleGAFA} applies to the newly constructed
expander on this manifold and supplies a potential $f$ for
\eqref{eq:soliton}. Replacing it by
\[
 f^{G}=\frac1{|G|}\sum_{\sigma\in G}f\circ S_\sigma
\]
preserves \eqref{eq:soliton}, since every $S_\sigma$ is an isometry.
Translate $f^G$ to its canonical potential $F$ as in
\eqref{eq:canonical-F}. It is invariant, proper, and strictly convex.

\emph{The Morse chart with the prescribed angular labeling.}
We explain the additional point in the coordinate construction of
Gianniotis--Schulze cited above. An arbitrary equivariant
$j:X\to\{F=R^2/4\}$ and \eqref{eq:morse-formula} would give a Morse
chart, but its limit link metric could be just a pullback of $\gamma$.
To recover the prescribed labeling, start instead on levels tending
to infinity, where the reference chart already specifies $\gamma$.

Their level-comparison construction gives equivariant normal-graph
identifications $K_a$ from $\{F_0=a\}$ to $\{F=a\}$ for large $a$.
Normal graphs are equivariant because their defining level functions and
the reference normal direction are invariant, hence the uniqueness of the graph
function makes it invariant. Flow the resulting parametrization to
other regular levels as
\begin{equation}\label{eq:high-level-chart}
 H_a(r,x)=
 \Fl_{\nabla F/|\nabla F|^2}^{\,r^2/4-a}
       \bigl(K_a(\Phi_0(2\sqrt a,x))\bigr).
\end{equation}
Each map is exactly equivariant, since both the graph identification and
the flow commute with the actions.
The rescaled metric and graph estimates in
\cite[equations~(6-1)--(6-6), pp.~3970--3973]{GianniotisSchulze}
give a convergent subsequence as the starting level tends to infinity.
Its limit is an asymptotic diffeomorphism $\Phi$ in these angular labels.
The construction compares with the specified cone at the starting
levels, and its limiting angular metric is $\gamma$, not an
unspecified relabeling of $\gamma$.
These estimates and constructions use smooth graphs and ordinary
differential equations on the covering manifolds. Normal exponential
maps and uniqueness of flows remain equivariant at points with
nontrivial stabilizers, so no freeness is needed.

Fixing a regular radius $R$ in the limiting chart and setting
$j(x)=\Phi(R,x)$, its expression is exactly
\[
 \Phi(r,x)=
 \Fl_{\nabla F/|\nabla F|^2}^{(r^2-R^2)/4}(j(x)).
\]
Indeed all the approximating maps \eqref{eq:high-level-chart} solve the
same flow equation. As computed in \eqref{eq:morse-formula},
$F\circ\Phi=r^2/4$ and $d\Phi(\partial_r)=
(r/2)\nabla F/|\nabla F|^2$. Equivariance passes to the limit.
Finally the sharp asymptotic estimates
\cite[Theorem~1.2 and proof of Theorem~3.2]{DeruellePisa}
give \eqref{eq:rate-two} in this Morse gauge.
\end{proof}

\subsection{Marked closedness}\label{subsection 3.2}

\begin{proposition}\label{prop:marked-closedness}
Let $X^{N-1}$ be closed and simply connected, $N\geq3$, and fix a
smooth finite-group action $\rho:G\to\Diff(X)$. Suppose
$h_i\to h_*$ smoothly, all these metrics are $G$-invariant, and
$\cR_{h_i}>I$, $\cR_{h_*}>I$.
If each $(X,h_i,\rho)$ admits an equivariant smoothing with a rate-two
Morse-flow chart as in \cref{thm:equivariant-smoothing}, then so does
$(X,h_*,\rho)$.
\end{proposition}

\begin{proof}
\emph{Step 1: compactness and the constants in the potentials.}
Denote by $(M_i,g_i,f_i,p_i)$ the expanders with the integral
normalization \eqref{eq:integral-normalization}.
The convergent link family is smoothly compact, and its curvature operators are uniformly strictly
above $I$ after discarding finitely many terms. It has uniform curvature-derivative bounds and a positive
volume lower bound. Bonnet--Myers gives
$\operatorname{diam}(X,h_i)\leq\pi$, and \eqref{eq:avr-link} gives
\[
 N\cdot\AVR(g_i)=\Vol_{h_i}(X)\geq v>0.
\]
The bounds explained after \cref{thm:deruelle-compactness} therefore
apply. A subsequence converges smoothly on pointed compact sets to
\[
 (M_\infty,g_\infty,f_\infty,p_\infty),\qquad
 \cR_{g_\infty}\geq0.
\]
Set
\begin{equation}\label{eq:Fi-sequence}
 c_i=\Scal_{g_i}(p_i)-f_i(p_i),\qquad F_i=f_i+c_i.
\end{equation}
Then $c_i\to c_\infty$ by pointed smooth convergence, and
$F_i\to F_\infty=f_\infty+c_\infty$ on the same compact sets.
The exact identity is
\[
 |\nabla F_i|^2+\Scal_{g_i}=F_i,
\]
including $i=\infty$. Adding $c_i$ fixes the radial normalization uniformly and removes a potential
constant from every transport equation below. It does not alter
the metrics, vector fields, or group actions.

For the limit, properness and
$\Hess F_\infty=\Ric_{g_\infty}+g_\infty/2\geq g_\infty/2$
give one nondegenerate minimum and
\begin{equation}\label{eq:limit-topology}
 M_\infty\cong\R^N
\end{equation}
as explained in \eqref{eq:expander-topology}. This conclusion only requires
nonnegative curvature.
We establish it now because simple connectivity is needed later.

\emph{Step 2: put the finite-annulus convergence in a common Morse gauge.}
Choose $R$ beyond the common radial threshold in
\eqref{eq:uniform-curvature}, and large enough that $R^2/4$ exceeds
all the minimum values $F_i(p_i)$. Let
\[
 X_0=\{F_\infty=R^2/4\},
\]
which is a closed connected regular level because $F_\infty$ is strictly convex and proper. In pointed convergence
coordinates, the corresponding $F_i$ levels are normal graphs over
$X_0$. Their graph parametrizations $j_i:X_0\to
\{F_i=R^2/4\}$ converge smoothly to the inclusion $j_\infty$.
All these levels lie in a common bounded pointed region by strict
convexity. Hence the pointed identifications cover the whole levels.

With $W_i=\nabla F_i/|\nabla F_i|^2$, define
\begin{equation}\label{eq:Ei-flow}
 E_i(r,x)=\Fl_{W_i}^{(r^2-R^2)/4}(j_i(x)),\qquad r\geq R.
\end{equation}
The usual smooth dependence of solutions of an ordinary differential equation on coefficients and initial data gives
\begin{equation}\label{eq:finite-annuli}
 E_i^*g_i\longrightarrow E_\infty^*g_\infty
 \quad\text{smoothly on }[R,L]\times X_0
 \quad\text{for every finite }L.
\end{equation}
The convergence in \eqref{eq:finite-annuli} is in the pointed identifications on the
relevant compact regions. Properness makes those regions compact.

The initially prescribed charts $\Phi_i$ may have different lower
cutoffs. First express each in its canonical $F_i$ radius. An additive
change of potential only changes $r^2$ by a constant, preserving the
link and the rate-two property. Then extend it inward by its Morse
flow to the common level $R^2/4$. This is possible because there are
no critical points above the minimum. Thus both $E_i$ and $\Phi_i$
can be compared on the common domain even if the original cutoffs
were not uniform.

\emph{Step 3: define $q_i(r)$ and prove a uniform tail estimate.}
For $E_{i,r}(x)=E_i(r,x)$, define the family of metrics on the
\emph{fixed} compact manifold $X_0$ by
\begin{equation}\label{eq:qi-definition}
 q_i(r)=r^{-2}E_{i,r}^*g_i,
\end{equation}
where the pullback is the induced metric on the level $F_i=r^2/4$, and the
division by $r^2$ removes the leading quadratic growth of a cone.
Thus $q_i(r)$ records the angular geometry of that level at unit scale.
%It is a tensor-valued function of $r$, not a point or an additional map between the manifolds.

Write
\[
 S_i(r,x)=\Scal_{g_i}(E_i(r,x)),\qquad
 u_i(r,x)=\frac{r^2}{4}-S_i(r,x)=|\nabla F_i|^2(E_i(r,x)).
\]
The gradient-flow direction is normal to the level sets, so there
are no radial-angular cross terms and
\begin{equation}\label{eq:metric-in-level-gauge}
 E_i^*g_i=\frac{r^2}{4u_i}\,dr^2+r^2q_i(r).
\end{equation}
On tangential vectors the derivative of the induced metric is
$(r/u_i)E_{i,r}^*\Hess F_i$. Inserting the soliton equation and
differentiating $r^{-2}$ gives the exact transport identity
\begin{align}
 \partial_rq_i
 &=\left(-\frac2r+\frac r{2u_i}\right)q_i
       +\frac1{ru_i}E_{i,r}^*\Ric_{g_i}\notag\\
 &=\frac{8S_i}{r(r^2-4S_i)}q_i
       +\frac4{r(r^2-4S_i)}E_{i,r}^*\Ric_{g_i}.
 \label{eq:qi-transport}
\end{align}
The uniform curvature estimates give $S_i=O(r^{-2})$.
The first term in the second line has $q_i$-norm $O(r^{-5})$.
The tangential Ricci tensor has $q_i$-norm $O(1)$, since its ambient
norm is $O(r^{-2})$ and the angular metric has been divided by $r^2$.
The second term therefore has $q_i$-norm $O(r^{-3})$.

We spell out why the differentiated estimates have constants independent
of $i$ and can be used in one fixed norm,  using the same argument as in the proof of
\cite[Theorem~3.2, pp.~508--509]{DeruellePisa}. Let $g_i^T$ be the induced metric on $\Sigma_{i,r}:=\{F_i=r^2/4\}$.
The shape endomorphism of this level set satisfies
\[
 A_i
 =\frac{(\Ric_{g_i}|_{T\Sigma_{i,r}\times T\Sigma_{i,r}})
                 ^{\sharp_{g_i^T}}+\frac12\id}{|\nabla F_i|}
 =r^{-1}\id+O(r^{-3}).
\]
Combine this with \eqref{eq:uniform-curvature} and the Gauss formulas.
Tangential derivatives of the tensors in \eqref{eq:qi-transport},
measured in the unit-scale metric $q_i(r)$, are uniformly bounded
at every order, with the displayed powers of $r$ for their
coefficients. Differentiating \eqref{eq:qi-transport} therefore gives
\[
 |(\nabla^{q_i(r)})^k\partial_rq_i|_{q_i(r)}
 \leq C_kr^{-3}\qquad(k\geq0).
\]
%These are the differentiated transport estimates used in the proof of \cite[Theorem~3.2, pp.~508--509]{DeruellePisa}; their uniform input is \cite[Proposition~4.3]{DeruellePisa}.

Fix a background metric $q_0$ on $X_0$. Since $q_i(R)$ converge
smoothly, their $C^k(q_0)$ norms are uniformly controlled initially.
The $k=0$ estimate and integration of $r^{-3}$ give uniform metric
equivalence of $q_i(r)$ and $q_0$ for all $r\geq R$.
For the connections, the variation formula is
\[
 \partial_r\Gamma(q_i)^\alpha_{\beta\gamma}
 =\frac12q_i^{\alpha\delta}
 \bigl(\nabla_\beta\partial_rq_{i,\gamma\delta}
      +\nabla_\gamma\partial_rq_{i,\beta\delta}
      -\nabla_\delta\partial_rq_{i,\beta\gamma}\bigr),
\]
where $\nabla=\nabla^{q_i(r)}$.
Integration controls the difference from the connection at $R$.
Differentiating this identity and proceeding inductively controls all
higher connection differences and metric derivatives. At each stage
the new coefficient is integrable in $r$, and the integration, or equivalently
Gronwall's inequality, gives a uniform bound.
Thus the moving-metric estimates imply the fixed-background estimates
\begin{equation}\label{eq:uniform-qi}
 \|\partial_rq_i(r)\|_{C^k(X_0,q_0)}
 \leq C_kr^{-3}\qquad(r\geq R).
\end{equation}
In particular, we do not assume convergence of the link metrics in these
coordinates in order to justify the norms.

For $t\geq r\geq R$, integration now gives
\begin{equation}\label{eq:qi-cauchy}
 \|q_i(t)-q_i(r)\|_{C^k(q_0)}
 \leq \frac{C_k}{2}(r^{-2}-t^{-2}).
\end{equation}
This proves uniform continuity in $r$, uniformly over the sequence,
and the stronger uniform Cauchy property as $r\to\infty$.
Consequently there are smooth positive metrics $b_i$ on $X_0$ with
\begin{equation}\label{eq:uniform-tail}
 b_i=\lim_{r\to\infty}q_i(r),\qquad
 \|q_i(r)-b_i\|_{C^k(q_0)}\leq\frac{C_k}{2r^2}.
\end{equation}
Positivity follows from the uniform two-sided metric equivalence.
The same estimates hold for $i=\infty$ by
\eqref{eq:finite-annuli}.

We can now interchange the sequence limit with the radial limit.
For every fixed $r$, $q_i(r)\to q_\infty(r)$ smoothly by
\eqref{eq:finite-annuli}, whereas
\[
 \|b_i-b_\infty\|_{C^k(q_0)}
 \leq \frac{C_k}{r^2}
       +\|q_i(r)-q_\infty(r)\|_{C^k(q_0)}.
\]
Choose $r$ large first and then $i$ large. It follows that
\begin{equation}\label{eq:bi-convergence}
 b_i\longrightarrow b_\infty\quad\text{smoothly}.
\end{equation}
The smooth convergence on finite annuli alone, or the uniform continuity
alone, would not justify this step. The uniform tail bound \eqref{eq:uniform-tail} is essential.
Together with \eqref{eq:metric-in-level-gauge} and its differentiated
estimates, it also identifies $(X_0,b_\infty)$ as the smooth rate-two
link in the chart $E_\infty$.

\emph{Step 4: recover the prescribed marking and pass the action to the end.}
Both $E_i$ and $\Phi_i$ now use the same canonical $F_i$ and solve the
same Morse equation. By \eqref{eq:morse-transition}, there is a
diffeomorphism $\psi_i:X_0\to X$, independent of $r$, such that
\begin{equation}\label{eq:morse-end-convergence}
 E_i=\Phi_i\circ(\id\times\psi_i),\qquad b_i=\psi_i^*h_i.
\end{equation}
The second identity follows by taking the rescaled tangential limit.
Apply \cref{lem:relabel} using \eqref{eq:bi-convergence} and
$h_i\to h_*$. Passing to a subsequence gives, smoothly,
\begin{equation}\label{eq:relabel-limit}
 \psi_i\to\psi_\infty,\qquad
 \psi_i^{-1}\to\psi_\infty^{-1},\qquad
 \psi_\infty^*h_*=b_\infty.
\end{equation}
In the $E_i$ coordinates put
\[
 \rho_i(\sigma)=\psi_i^{-1}\rho(\sigma)\psi_i,\qquad
 \rho_\infty(\sigma)=\psi_\infty^{-1}\rho(\sigma)\psi_\infty.
\]
Exact equivariance of the prescribed chart and the global isometry
property imply
\[
 (\id\times\rho_i(\sigma))^*E_i^*g_i=E_i^*g_i.
\]
Using finite-annulus convergence and smooth convergence of $\rho_i$,
we may pass to the limit to obtain
\[
 (\id\times\rho_\infty(\sigma))^*E_\infty^*g_\infty
 =E_\infty^*g_\infty.
\]
This defines an isometric action on the connected end
$E=E_\infty((R,\infty)\times X_0)$.
%No quotient has been taken, and no freeness has been used.

\emph{Step 5: strict positivity, analyticity, and global extension.}
By \eqref{eq:relabel-limit}, the limiting link satisfies
$\cR_{b_\infty}>I$.
The limiting manifold is simply connected by
\eqref{eq:limit-topology}. Therefore the strong-maximum-principle and
holonomy argument in
\cite[proof of Theorem~1.3, pp.~244--245]{DeruelleGAFA}
applies to this nonnegatively curved expander and its strictly curved
link and yields $\cR_{g_\infty}>0$.

To justify the analytic continuation step, apply Ivey's theorem
\cite[Theorem~1]{Ivey} to the \emph{limiting} equation
\[
 \mathcal L_{\nabla F_\infty}g_\infty
 =2\Ric_{g_\infty}+g_\infty.
\]
It follows that $g_\infty$ is real analytic in harmonic coordinates.
%We do not infer analyticity from the smooth convergence of $g_i$: such an inference is false for general sequences of analytic metrics. The soliton equation is the reason analyticity holds here.

For $\sigma\in G$, define the end isometry
\[
 a_\sigma=E_\infty\circ(\id\times\rho_\infty(\sigma))
                      \circ E_\infty^{-1}:E\to E.
\]
The manifold $M_\infty$ is complete and simply connected.  %, connected,
Both fundamental-group surjectivity requirements of
\cite[Theorem~5.2]{KotschwarWang} hold because their target group is
$\pi_1(M_\infty)=0$. Hence \cref{lem:globalize-action} extends
$a_\sigma$ to a unique global isometry
$\widehat\rho_\infty(\sigma)$.
The compositions $\widehat\rho_\infty(\sigma)
\widehat\rho_\infty(\tau)$ and
$\widehat\rho_\infty(\sigma\tau)$ agree on $E$, and therefore everywhere. Moreover,
the identity and inverse laws follow as well. These extensions thus
form a genuine action.

On the end $F_\infty=r^2/4$ and $a_\sigma$ preserves $r$.
Moreover, the global soliton equation gives
\[
 \Hess_{g_\infty}
 \bigl(F_\infty\circ\widehat\rho_\infty(\sigma)-F_\infty\bigr)=0.
\]
Its differential is parallel and zero on $E$, so the difference
vanishes globally. Thus the potential is invariant and its unique
minimum is fixed.

\emph{Step 6: restore the original action, not just a conjugate.}
Set
\[
 \widetilde\Phi_\infty
 =E_\infty\circ(\id\times\psi_\infty^{-1}):
 (R,\infty)\times X\longrightarrow E.
\]
Its link metric is $(\psi_\infty^{-1})^*b_\infty=h_*$.
The identity $\rho_\infty(\sigma)\psi_\infty^{-1}
=\psi_\infty^{-1}\rho(\sigma)$ gives
\[
 \widehat\rho_\infty(\sigma)\circ\widetilde\Phi_\infty
 =\widetilde\Phi_\infty\circ(\id\times\rho(\sigma)).
\]
The relabeling is independent of $r$, so it preserves
$F_\infty\circ\widetilde\Phi_\infty=r^2/4$, the Morse-flow equation,
and all rate-two estimates. This is the required marked smoothing.
\end{proof}

\subsection{An invariant path and completion of the proof}

\begin{lemma}\label{lem:invariant-path}
Let $(X^m,h)$ be closed and simply connected, $m\geq2$, with
$\cR_h>I$, and let $G$ act isometrically on $(X,h)$.
There is a smooth $G$-invariant path $(h_s)_{s\in[0,1]}$ with
\[
 h_0=h,\qquad h_1=c^2g_{\rm rd},\qquad
 \cR_{h_s}>I,\qquad 0<c<1,
\]
after an identification of $X$ with $\Sph^m$.
Its image is compact in the smooth topology.
\end{lemma}

\begin{proof}
Run the volume-normalized Ricci flow
\begin{equation}\label{eq:normalized-flow}
 \partial_t h(t)=-2\Ric_{h(t)}
      +\frac2m\overline{\Scal}(h(t))\,h(t),\qquad h(0)=h,
\end{equation}
where
\[
 \overline{\Scal}(k)=
 \frac{\int_X\Scal_k\,d\mu_k}{\Vol_k(X)}
\]
is the average scalar curvature.
For $m\geq3$, B\"ohm--Wilking \cite{BohmWilking} prove preservation of
positive curvature operator and smooth convergence to a
constant-positive-curvature metric. For $m=2$, use Hamilton's
positive-curvature surface theorem \cite{HamiltonSurface}.
Simple connectivity makes the limit a round sphere, not a nontrivial
spherical space form.

Here is the precise meaning of \emph{naturality and uniqueness}.
For each fixed diffeomorphism $\varphi=\rho(\sigma)$,
\[
 \Ric_{\varphi^*k}=\varphi^*\Ric_k,\qquad
 \overline{\Scal}(\varphi^*k)=\overline{\Scal}(k).
\]
Thus $\varphi^*h(t)$ satisfies \eqref{eq:normalized-flow} and has the
same initial metric as $h(t)$. 
\cite[Theorem~1]{KotschwarUniqueness} proves uniqueness for complete bounded-curvature Ricci flows and recalls Hamilton's compact case. Such uniqueness gives $\varphi^*h(t)=h(t)$. Actually, 
%Uniqueness gives $\varphi^*h(t)=h(t)$.
%For a precise uniqueness reference, use \cite[Theorem~1]{KotschwarUniqueness}, which proves uniqueness for complete bounded-curvature Ricci flows and recalls Hamilton's compact case. 
on a closed manifold curvature is bounded on each finite time
interval. The normalized statement follows by the usual scalar metric
rescaling and time change to the unnormalized Ricci flow.
For the two flows just compared, the average scalar curvatures agree
by pullback invariance, so these changes are the same.
This proves invariance for every finite $t$, and smooth convergence
proves invariance of the limit $h_\infty$.

We need a smooth path on a \emph{finite} parameter interval.
Choose a large finite time $T$ such that $h(T)$ lies in a small
$C^2$-neighborhood of $h_\infty$. That neighborhood can be chosen
so that the straight segment
\[
 (1-u)h(T)+u h_\infty,\qquad 0\leq u\leq1,
\]
consists of positive-curvature-operator metrics.
Both endpoints are invariant, so the segment is invariant.
Concatenate the flow from $0$ to $T$ with this segment and use smooth
parameter functions flat at the joining point. This yields a smooth
invariant path $\bar h_s$ from $h$ to $h_\infty$ with
$\cR_{\bar h_s}>0$. 
%No smooth compactification of infinite flow time is needed.

The normalized flow preserves volume. Also, compactness and
$\cR_h>I$ give $\varepsilon>0$ with
$\Ric_h\geq(m-1)(1+\varepsilon)h$. Bishop comparison gives
\[
 \Vol_h(X)\leq(1+\varepsilon)^{-m/2}
                      \Vol(\Sph^m,g_{\rm rd})
 <\Vol(\Sph^m,g_{\rm rd}).
\]
Therefore $h_\infty=c^2g_{\rm rd}$ with $c<1$.
Let $\lambda(s)$ be the minimum, over $X$, of the least eigenvalue of
$\cR_{\bar h_s}$. It is positive and continuous, and greater than one
near both endpoints. Choose a smooth positive $b(s)$ equal to one
near the endpoints and satisfying $b(s)^2<\lambda(s)$ everywhere.
Then
\[
 h_s=b(s)^2\bar h_s,\qquad
 \cR_{h_s}=b(s)^{-2}\cR_{\bar h_s}>I.
\]
This path has the required endpoints and symmetry.
Its smooth image over a compact interval is a compact smooth family.
\end{proof}

\begin{proof}[Proof of \cref{thm:equivariant-smoothing}]
Choose the path in \cref{lem:invariant-path}. Let
$\mathcal S\subset[0,1]$ consist of parameters $s$ for which
$(X,h_s,\rho)$ admits an equivariant smoothing with a rate-two
Morse-flow chart and the original marking.

The round endpoint belongs to $\mathcal S$.
Choose an isometry
$J:(X,h_1)\to(\Sph^m,c^2g_{\rm rd})$.
The conjugated action $J\rho J^{-1}$ is orthogonal.
The rotationally symmetric Bryant expander asymptotic to
$dr^2+c^2r^2g_{\rm rd}$ has positive curvature operator when $c<1$;
see \cite{DeruelleGAFA}. Its warped-product metric and radial
potential are invariant under every orthogonal transformation.
If $t$ is its distance radius and $F(t)$ its canonical potential, use
$t=t(r)$ determined by $F(t(r))=r^2/4$ on the end.
Differentiation gives
\[
 \frac{dt}{dr}=\frac{r}{2F'(t)},\qquad
 d\Phi_B(\partial_r)
 =\frac r2\frac{\nabla F}{|\nabla F|^2}\circ\Phi_B,
 \quad \Phi_B(r,\theta)=(t(r),\theta).
\]
The chart $\Phi_B\circ(\id\times J)$ has the original action $\rho$
as its marking on $X$.

By \cref{prop:open}, $\mathcal S$ is relatively open.
By \cref{prop:marked-closedness}, it is relatively closed, since the
path is a compact smooth family inside $\{\cR>I\}$.
It is nonempty and $[0,1]$ is connected, hence
$\mathcal S=[0,1]$. The member at $s=0$ proves the theorem.
Rename its canonically translated potential $f$ to obtain the notation
of the statement.
\end{proof}

\section{Fixed components and the volume comparison}\label{sec:fixed}

We now apply \cref{thm:equivariant-smoothing} to prove
\cref{thm:fixed-link,thm:main}.  The geometric point is that an exactly
equivariant asymptotic chart identifies each fixed component of the link
with the link of a complete totally geodesic submanifold.  One such conical
end already gives the lower volume-growth bound needed for comparison;
there is no need to exclude other ends of the fixed submanifold.

\subsection{Fixed submanifolds and their ends}

\begin{lemma}\label{lem:fixed-submanifold}
Let $(M,g)$ be a Riemannian manifold and let
$\mathcal A\subset\Isom(M,g)$.  Every connected component $Y$ of the
common fixed set
\[
 M^{\mathcal A}
 =\{p\in M:\sigma(p)=p\text{ for every }\sigma\in\mathcal A\}
\]
is a closed embedded totally geodesic submanifold, with
\[
 T_pY=\{v\in T_pM:d\sigma_p(v)=v
                 \text{ for every }\sigma\in\mathcal A\}.
\]
If $(M,g)$ is complete, then $Y$ is complete for its induced metric.
\end{lemma}

\begin{proof}
This is the classical fixed-point theorem for isometries; see \cite[Chapter~II, Theorem~5.1, p.~59]{Kobayashi}.
This theorem applies to an arbitrary set of isometries.
We recall the local argument to identify the tangent spaces
and to establish completeness when the ambient manifold is complete.
We recall the local argument because it also explains the completeness assertion.  At a common
fixed point $p$, the exponential map satisfies
\[
 \sigma\circ\exp_p=\exp_p\circ d\sigma_p.
\]
On a sufficiently small normal ball, $\exp_p$ is injective.  Its common
fixed points are therefore precisely the exponential image of the vector
subspace fixed by all $d\sigma_p$.  This proves embeddedness and the
tangent-space formula.  The common fixed set is closed, and each of its
connected components is closed as well.

If $v\in T_pY$, the geodesic with initial data $(p,v)$ is fixed by every
$\sigma\in\mathcal A$, by uniqueness of geodesics.  Its connected image
stays in $Y$, hence is totally geodesic.  When $M$ is complete, this geodesic
exists for all time.  Thus the induced metric on $Y$ is geodesically
complete, and hence complete by Hopf--Rinow.
\end{proof}

\begin{lemma}\label{lem:fixed-end}
Let $\rho:G\to\Isom(X,h)$ and
$\widehat\rho:G\to\Isom(M,g)$ be finite-group actions, where $(M,g)$
is complete.  Suppose a smooth asymptotic chart
\[
 \Phi:(R_0,\infty)\times X\longrightarrow M\setminus K
\]
satisfies
\[
 \widehat\rho(\sigma)\circ\Phi
 =\Phi\circ(\id\times\rho(\sigma)),\qquad \sigma\in G.
\]
If $Y^d$ is a connected component of $X^G$, then there is a unique
connected component $\Sigma^{d+1}$ of $M^G$ containing
$\Phi((R_0,\infty)\times Y)$.  It is complete and totally geodesic.
For sufficiently large $R_1$, the set
\[
 E_Y=\Phi((R_1,\infty)\times Y)
\]
is an end of $\Sigma$, smoothly asymptotic to the cone with metric
$dr^2+r^2h|_Y$.  The component $\Sigma$ may have other ends.
\end{lemma}

\begin{proof}
The equivariance and injectivity give the identity
\begin{equation}\label{eq:fixed-end-identity}
 M^G\cap(M\setminus K)
 =\Phi((R_0,\infty)\times X^G).
\end{equation}
Indeed, $(r,x)$ maps to a common fixed point if and only if every
$\rho(\sigma)$ fixes $x$.  The connected set
$\Phi((R_0,\infty)\times Y)$ lies in a unique component $\Sigma$ of
$M^G$.  The local identity \eqref{eq:fixed-end-identity} gives
$\dim\Sigma=d+1$. By \cref{lem:fixed-submanifold}, $\Sigma$ is complete and totally geodesic. 

Truncate the asymptotic chart at $R_1>R_0$ and let
\[
 K_1=M\setminus\Phi((R_1,\infty)\times X),
\]
which is $G$-invariant. 
As in the exterior-chart convention of \cref{sec:inputs}, $K_1$ is also
compact. The set $E_Y$ is a connected
component of $\Sigma\setminus K_1$, and it is diffeomorphic to
$(R_1,\infty)\times Y$.  Its tails, obtained by increasing $R_1$,
represent an end because $Y$ is compact and connected.  
%No soliton potential is needed for this truncation.

%The expression \emph{fixed cylinder} refers to this product parametrization: $E_Y$ has a radial coordinate and a cross-section diffeomorphic to $Y$.  It does not assert that its induced metric is the product metric $dr^2+h|_Y$.  Even the limiting cone metric is the warped product $dr^2+r^2h|_Y$.

The fixed cylinder $E_Y$ has a radial coordinate and a cross-section diffeomorphic to $Y$. But its metric is not the product metric $dr^2+h|_Y$. Even the limiting cone metric is the warped product $dr^2+r^2h|_Y$. 
To check the metric assertion, we write
\[
 g_C=dr^2+r^2h,\qquad e=\Phi^*g-g_C,
 \qquad g_{C(Y)}=dr^2+r^2h|_Y.
\]
Let $\iota:(R_1,\infty)\times Y\hookrightarrow
(R_1,\infty)\times X$ denote the inclusion and set
$\Phi_Y=\Phi\circ\iota$.  Then
\begin{equation}\label{eq:fixed-end-error}
 \Phi_Y^*(g|_\Sigma)=g_{C(Y)}+e_Y,
 \qquad e_Y=\iota^*e.
\end{equation}
%The phrase \emph{relative norm} means the tensor norm measured with the cone metric. In particular,
It can be deduced from
\[
 |e_Y|_{g_{C(Y)}}\leq |e|_{g_C}\circ\iota=o(1)
\]
that for every $\varepsilon>0$, on a sufficiently distant tail,
\[
 (1-\varepsilon)g_{C(Y)}
 \leq\Phi_Y^*(g|_\Sigma)
 \leq(1+\varepsilon)g_{C(Y)}.
\]
It is an asymptotic metric error estimate in the relative norm, which does not mean that all angular coordinate coefficients of $e_Y$ tend to zero, since angular vectors have cone lengths proportional to $r$.
%This is a relative error statement: angular vectors have cone lengths proportional to $r$, so small relative error does not mean that all angular coordinate coefficients of $e_Y$ tend to zero.

%The \emph{differentiated estimates} are the corresponding bounds on every cone covariant derivative of $e_Y$.  
The fixed component $Y$ is totally
geodesic in $X$, so $\iota$ is totally geodesic for the cone metrics.
Consequently, tangential restriction commutes with covariant
differentiation, and for every integer $j\geq0$,
\[
 |(\nabla^{g_{C(Y)}})^je_Y|_{g_{C(Y)}}
 \leq |(\nabla^{g_C})^je|_{g_C}\circ\iota.
\]
The asymptotic estimates in \cref{sec:inputs} therefore restrict to this
end.  In particular, a rate-two ambient estimate gives
$|(\nabla^{g_{C(Y)}})^je_Y|_{g_{C(Y)}}\leq C_jr^{-2-j}$.
This proves the asserted smooth conical asymptotics.
\end{proof}

\subsection{Volume growth from one conical end}

\begin{lemma}\label{lem:end-volume}%\label{lem:volume-growth}
Let $(\Sigma^k,g_\Sigma)$ be complete and connected.  If one of its ends
is smoothly asymptotic to the cone over a closed connected manifold
$(Y^{k-1},h_Y)$, then for every $p\in\Sigma$,
\begin{equation}\label{eq:end-volume}
 \liminf_{R\to\infty}
 \frac{\Vol B_\Sigma(p,R)}{R^k}
 \geq\frac{\Vol_{h_Y}(Y)}{k}.
\end{equation}
Other ends of $\Sigma$ are allowed.  
%Only the zeroth-order relative metric estimate is needed for this conclusion.
\end{lemma}

\begin{proof}
Let $\Phi$ parametrize the chosen end, write $g_0=dr^2+r^2h_Y$, and set
$\widehat g=\Phi^*g_\Sigma$.  Given $\delta\in(0,1)$, choose
$R_\delta$ large enough that
\[
 (1-\delta)g_0\leq\widehat g\leq(1+\delta)g_0
 \qquad(r\geq R_\delta).
\]
Let $j_\delta:Y\to\Sigma$ parametrize the cross-section at
$r=R_\delta$.  Its image is compact, so
\[
 C_\delta=\max_{x\in Y}d_\Sigma(p,j_\delta(x))<\infty.
\]
Joining $p$ to $j_\delta(x)$ and then following the radial curve at
fixed $x$ gives
\[
 d_\Sigma(p,\Phi(r,x))
 \leq C_\delta+\sqrt{1+\delta}\,(r-R_\delta).
\]
Hence, for $R>C_\delta$, the image of the whole annular region
\[
 R_\delta\leq r<
 R_\delta+\frac{R-C_\delta}{\sqrt{1+\delta}}
\]
lies in $B_\Sigma(p,R)$.

The lower metric bound implies the lower volume-density bound
\[
 d\mu_{\widehat g}
 \geq(1-\delta)^{k/2}r^{k-1}\,dr\,d\mu_{h_Y}.
\]
Integrating over the annular region, dividing by $R^k$, and sending
$R\to\infty$ yields
\[
 \liminf_{R\to\infty}\frac{\Vol B_\Sigma(p,R)}{R^k}
 \geq\frac{\Vol_{h_Y}(Y)}{k}
       \left(\frac{1-\delta}{1+\delta}\right)^{k/2}.
\]
Letting $\delta\downarrow0$ proves \eqref{eq:end-volume}. The argument
uses only points in the chosen end, and any other ends can only increase the
volume of the ball.
\end{proof}

For clarity, when $Y$ is a circle of length $L$, arclength gives
$Y=\R/L\mathbb Z$ with metric $ds^2$.  Its cone is two-dimensional with metric and volume form as
\[
 g_{C(Y)}=dr^2+r^2ds^2,\qquad d\mu_{C(Y)}=r\,dr\,ds.
\]
In particular,
\[
 \operatorname{Area}\{R_0<r<R\}
 =\int_{R_0}^R\int_0^L r\,ds\,dr
 =\frac L2(R^2-R_0^2).
\]
Thus the circle length contributes the coefficient $L/2$ to quadratic
area growth.  This agrees with the unnormalized convention for $\AVR$
in \cref{sec:inputs}. With normalization by the Euclidean disk area
$\pi R^2$, the coefficient would instead be $L/(2\pi)$.

\subsection{Proofs of the comparison theorems}

\begin{proof}[Proof of \cref{thm:fixed-link}]
First suppose $\cR_h>I$.  Apply \cref{thm:equivariant-smoothing}, and
let $\Sigma^{d+1}$ be the complete totally geodesic fixed component
provided by \cref{lem:fixed-end}.  The ambient curvature operator is
positive.  Since the second fundamental form of $\Sigma$ vanishes, the
Gauss equation identifies its sectional curvatures with the corresponding
ambient sectional curvatures.  In particular, $\Ric_\Sigma\geq0$.
\Cref{lem:end-volume} gives
\[
 \AVR(g|_\Sigma)
 \geq\frac{\Vol_h(Y)}{d+1}.
\]
Here the asymptotic volume ratio exists by Bishop--Gromov comparison,
which also gives, for every $R>0$,
\[
 \frac{\Vol B_\Sigma(p,R)}{R^{d+1}}
 \leq\omega_{d+1}
 =\frac{\Vol(\Sph^d,g_{\rm rd})}{d+1}.
\]
Combining the two inequalities proves
$\Vol_h(Y)\leq\Vol(\Sph^d,g_{\rm rd})$ in the strict case.

If $\cR_h\geq I$, let $h_c=c^2h$ for $0<c<1$.  Constant scaling
multiplies curvature-operator eigenvalues by $c^{-2}$, so
$\cR_{h_c}\geq c^{-2}I>I$.  The action is still isometric and its fixed
set is unchanged.  Applying the strict case to $h_c$ yields
\[
 c^d\Vol_h(Y)=\Vol_{h_c}(Y)
 \leq\Vol(\Sph^d,g_{\rm rd}).
\]
Letting $c\uparrow1$ proves the claimed weak-inequality case.
\end{proof}

\begin{proof}[Proof of \cref{thm:main}]
By \cref{lem:fixed-submanifold}, a one-dimensional fixed component is an embedded circle.  The totally geodesic property makes it a closed geodesic with its
induced length.  In dimension three every bivector is decomposable, so
$\sec_g\geq1$ is equivalent to $\cR_g\geq I$.  Apply \cref{thm:fixed-link} with
$G=\langle\tau\rangle$ and $d=1$ to obtain
\[
 \Length_g(Y)\leq\Vol(\Sph^1,g_{\rm rd})=2\pi.
\]
On the unit round three-sphere, the linear involution
\[
 (x_1,x_2,x_3,x_4)\longmapsto(x_1,x_2,-x_3,-x_4)
\]
fixes a great circle of length $2\pi$, proving sharpness.
\end{proof}

\section{Acknowledgements}
This paper is partially supported by the Jiangsu Provincial Scientific Research Center of Applied Mathematics (Grant No. BK20233002), NSFC (Grant Nos.  12301032, 12271093), NSF of Jiangsu (Grant No. BK20230803),  Fundamental Research Funds for the Central Universities (Grant No. 4007012402), and Zhishan Scholars Programs of Southeast University (Grant No. 2242024RCB0039).

\section*{Use of generative AI.} 
During the preparation of this work, the authors used ChatGPT  for exploratory computations, possible proof directions, and editorial polishing.  The authors reviewed and verified all output.

\end{document}